\magnification=\magstep1
\parindent=0pt
\baselineskip=13pt

\font\BF=cmbx12 at 17pt
\font\bfsl=cmbxsl9
\font\slsl=cmsl9
\font\am=msam10
\font\bm=msbm10

\def\eop{\hfill{\am\char"03}}
\def\loe{\mathrel{\hbox{\am\char"36}}}
\def\RN{\hbox{\bm\char"52}}
\setbox9=\hbox{\rm(A) \ }
\def\CA{$\hbox to\wd9{\rm(A) \ }$}
\def\CL{$\hbox to\wd9{\rm(\hfil\hfil L\hfil) \ }$}
\def\CR{$\hbox to\wd9{\rm(\hfil\hfil R\hfil) \ }$}
\def\Cr{$\hbox to\wd9{\rm(r)\hfil}$}


\centerline{\BF The last proof of extreme value theorem}
\medskip
\centerline{\BF and intermediate value theorem}

\bigskip\medskip


\hrule
\medskip
\centerline{Claude-Alain Faure}
\medskip
\hrule

\bigskip\medskip


{\bfsl Abstract.\ \ \slsl The paper gives a unified and simple proof of both theorems
and Cousin's theorem.}

\bigskip


{\bf Key-lemma.} \ Let $\cal I$ be a set of subintervals of a closed interval $[a,b]$.
We suppose that the set $\cal I$ satisfies the following conditions:
\smallskip
(A) \ $[r,s]\in\cal I$ and $[s,t]\in{\cal I}\ \Rightarrow\ [r,t]\in\cal I$\quad(additivity)
\smallskip
\CL for every $s\in{]a,b]}$, there exists $0<\delta\loe s-a$ such that $r\in{[s-\delta,s[}$
$\Rightarrow$ $[r,s]\in\cal I$
\smallskip
\CR for every $s\in{[a,b[}$, there exists $0<\delta\loe b-s$ such that $t\in{]s,s+\delta]}$
$\Rightarrow$ $[s,t]\in\cal I$
\smallskip
or the weaker condition
\smallskip
\Cr for every $s\in{[a,b[}$, there exists $t\in{]s,b]}$ with $[s,t]\in\cal I$
\smallskip
Then the whole interval $[a,b]$ must be an element of $\cal I$.

\medskip


{\bf Proof.} \ We consider the set $E=\{a\}\cup\{x\in{]a,b]}\mathrel/[a,x]\in{\cal I}\}$
and $s=\sup(E)$.
\smallskip
We first show that $s\in E$. We may assume that $s>a$. Let $0<\delta\loe s-a$
be as in (L). By definition of the supremum, there exists $r\in{]s-\delta,s]}\cap E$.
If $r=s$, we are done. And if $r<s$, then $[a,r]\in\cal I$ and $[r,s]\in{\cal I}$
$\Rightarrow$ $[a,s]\in\cal I$ by additivity. So $s\in E$.
\smallskip
We now show that $s=b$, which proves the lemma. Suppose on the contrary that
$s<b$. By (r) we have $[s,t]\in\cal I$ for some $t\in{]s,b]}$, and hence $[a,t]\in\cal I$
by additivity. So $t\in E$, in contradiction to the supremacy of $s$.\eop

\bigskip


{\bf Extreme value theorem.} \ Let $f:[a,b]\to\RN$ be a continuous function. Then
there exist $c,d\in[a,b]$ such that $f(c)\loe f(x)\loe f(d)$ for all $x\in[a,b]$.

\medskip


{\bf Proof.} \ Suppose on the contrary that such a number $d$ does not exist. We
show that
\medskip
\centerline{${\cal I}=\{I\subset[a,b]\mathrel/\hbox{there exists }u\in[a,b]\hbox{ such
that }f(x)<f(u)\hbox{ for all }x\in I\}$}
\medskip
satisfies the conditions of the key-lemma.
\smallskip
\CA Let $[r,s]\in\cal I$ and $[s,t]\in\cal I$ . Then $f(x)<f(u)$ for all $x\in[r,s]$ and
$f(x)<f(v)$ for all $x\in[s,t]\ \Rightarrow\ f(x)<\max{\{f(u),f(v)\}}$ for all $x\in[r,t]$.
So $[r,t]\in\cal I$.
\smallskip
\CL Let $s\in{]a,b]}$ . By assumption we have $f(s)<f(u)$ for some $u\in[a,b]$.
Choose $\varepsilon>0$ with $f(s)+\varepsilon<f(u)$. By continuity, there exists
$0<\delta\loe s-a$ such that $x\in[s-\delta,s]$ $\Rightarrow$ $f(x)\loe f(s)+
\varepsilon<f(u)$, and this shows that $[r,s]\in\cal I$ for all $r\in{[s-\delta,s[}$.
\smallskip
\CR is proved in the same way.
\smallskip
By the key-lemma we conclude that $[a,b]\in\cal I$, which is clearly impossible.
\hfill\eop

\bigskip


{\bf Intermediate value theorem.} \ Let $f:[a,b]\to\RN$ be a continuous function.
Then for every value $y$ between $f(a)$ and $f(b)$, there exists $c\in[a,b]$ with
$f(c)=y$.

\medskip


{\bf Proof.} \ Suppose on the contrary that such a number $c$ does not exist. We
show that
\medskip
\centerline{${\cal I}=\{I\subset[a,b]\mathrel/f(x)<y\hbox{ for all }x\in I\hbox{ or }f(x)>y
\hbox{ for all }x\in I\}$}
\medskip
satisfies the conditions of the key-lemma.
\smallskip
\CA Let $[r,s]\in\cal I$ and $[s,t]\in\cal I$ . Since $s\in{[r,s]}\cap{[s,t]}$, the inequality
must be the same on both intervals, and hence $[r,t]\in\cal I$.
\smallskip
\CL Let  $s\in{]a,b]}$ . We have $f(s)\not=y$, say $f(s)<y$. Choose $\varepsilon>0$
with $f(s)+\varepsilon<y$. By continuity, there exists $0<\delta\loe s-a$ such that $x
\in[s-\delta,s]$ $\Rightarrow$ $f(x)\loe f(s)+\varepsilon<y$, and this shows that $[r,s]
\in\cal I$ for all $r\in{[s-\delta,s[}$.
\smallskip
\CR is proved in the same way.
\smallskip
By the key-lemma we conclude that $[a,b]\in\cal I$, which is clearly impossible.\eop

\bigskip


\setbox8=\hbox{$\bullet$ \ }

In his famous paper [3], R.\ A.\ Gordon proved these theorems using Cousin's
theorem, a proof which is also exposed in [1].  Other elegant proofs of the extreme
value theorem can be found in [2], [4]. As another application of the key-lemma, we
give a very short proof of Cousin's theorem. We first recall the following definitions:
\smallskip
$\bullet$ \ a {\sl gauge} on $[a,b]$ is a positive function $\delta:[a,b]\to\RN$
\smallskip
$\bullet$ \ a {\sl tagged partition} of $[a,b]$ is a finite sequence $a=a_0<a_1<
\dots<a_{n-1}<a_n=b$
\par
$\kern\wd8$together with numbers $x_i\in[a_{i-1},a_i]$ for all $i=1,\dots,n$
\smallskip
$\bullet$ \ given a gauge $\delta$ on $[a,b]$, we say that a tagged partition is
{\sl$\delta$-fine} if
\par
$\kern\wd8{[a_{i-1},a_i]}\subset{[x_i-\delta(x_i),x_i+\delta(x_i)]}$ for all $i=1,
\dots,n$
\smallskip
$\delta$-fine tagged partitions will be called $\delta$-fine partitions for short.

\bigskip


{\bf Cousin's theorem.} \ For every gauge $\delta$ on $[a,b]$, there exists a
$\delta$-fine partition of $[a,b]$.

\medskip


{\bf Proof.} \ We show that ${\cal I}=\{I\subset[a,b]\mathrel/\hbox{there exists a
$\delta$-fine partition of }I\}$ satisfies the conditions of the key-lemma.
\smallskip
(A) \ It is immediate: just rename the points of the second partition.
\smallskip
\CL Let $s\in{]a,b]}$ . We consider $\eta=\min{\{
\delta(s),s-a\}}$. Then for every $r\in{[s-\eta,s[}$, the interval $[r,s]$ itself together
with the tag $s$ is a $\delta$-fine partition of $[r,s]$.
\smallskip
\CR is proved in the same way.
\smallskip
By the key-lemma we conclude that $[a,b]\in\cal I$, which proves the theorem.\eop

\bigskip


\setbox7=\hbox{[0] \ }

{\bf References}
\smallskip
$\hbox to\wd7{[1] \ }$R.\ G.\ Bartle and D.\ R.\ Sherbert, {\sl Introduction to real
analysis 4th edition,} New York, ${}\kern\wd7$John Wiley \& Sons, 2011.
\smallskip
$\hbox to\wd7{[2] \ }$S.\ J.\ Ferguson, {\sl A one-sentence line-of-sight proof of the
extreme value theorem,} Amer. ${}\kern\wd7$Math.\ Monthly {\bf121} (2014), 331.
\smallskip
$\hbox to\wd7{[3] \ }$R.\ A.\ Gordon, {\sl The use of tagged partitions in elementary
real analysis,} Amer.\ Math. ${}\kern\wd7$Monthly {\bf105} (1998), 107--117, 886.
\smallskip
$\hbox to\wd7{[4] \ }$M.\ H.\ Protter and C.\ B.\ Morrey, {\sl A first course in real
analysis,} New York, Springer, ${}\kern\wd7$1977.

\bigskip\medskip


\centerline{\vbox{\hbox{Gymnase de la Cit\'e, place de la Cath\'edrale 1, 1014
Lausanne, Switzerland}
\hbox{Email address: {\tt claudealain.faure@eduvaud.ch}}}}

\end